\documentclass[10pt]{amsart}
\usepackage{graphicx}
\usepackage{amsmath}
\usepackage{amssymb}
\usepackage{amsthm}
\usepackage[utf8]{inputenc}  
\usepackage[T1]{fontenc} 
\usepackage[french, german, italian, english]{babel} 
\usepackage{url}
\renewcommand{\Re}{\operatorname{Re}}
\renewcommand{\Im}{\operatorname{Im}}
\usepackage[author-year]{amsrefs}

\begin{document}
\title{Long-term history and ephemeral configurations}
\author{Catherine Goldstein}\thanks{I would like to warmly thank S\'ebastien Gauthier, Fran{\c c}ois L\^e, Ir\`ene Passeron and Jim Ritter for their helpful suggestions and commentaries.}
\address{CNRS, Institut de math\'ematiques de Jussieu-Paris Rive Gauche, Sorbonne Universit\'e,  Univ. Paris Diderot, Case 247, 4 place Jussieu, 75252 Paris Cedex, France} \email{catherine.goldstein@imj-prg.fr}
\thanks{I would like to warmly thank S\'ebastien Gauthier, Fran{\c c}ois L\^e, Ir\`ene Passeron and Jim Ritter for their helpful suggestions and commentaries.}


\maketitle

\begin{abstract}
Mathematical concepts and results have often been given a long  history, stretching far back in time. Yet recent work in the history of mathematics has tended to focus on local topics, over a short term-scale, and on the study of ephemeral configurations of mathematicians, theorems or practices. The first part of the paper explains why this change has taken place: a renewed interest in the connections between mathematics and society, an increased attention to the variety of components and aspects of mathematical work, and a critical outlook on historiography itself. The problems of a long-term history are illustrated and tested using a number of episodes in the nineteenth-century history of Hermitian forms, and finally, some open questions are proposed.
\end{abstract}

``Mathematics is the art of giving the same name to different things,'' wrote Henri Poincar\'e at the very beginning of the twentieth century (\cite[31]{Poincare1908}). The sentence, to be found in a chapter entitled ``The future of mathematics'' seemed particularly relevant around 1900: a structural point of view and a wish to clarify and to firmly found mathematics were then gaining ground and both contributed to shorten chains of argument and gather together under the same word phenomena which had until then been scattered (\cite{Corry2004}). Significantly, Poincar\'e's examples included uniform convergence and the concept of group.

\section{Long-term histories}
But the view of mathematics encapsulated by this --- that it deals somehow with\linebreak ``sameness'' --- has also found its way into the history of mathematics. It has been in particular a key feature (though often only implicitly) in the writing of most long-term histories. In the popular genres of the history of $\pi$ or of the Pythagorean theorem from Antiquity to present times, is hidden the idea that, despite changes in symbolism, despite the use or not of figures, tables or letters, despite the presence or not of proofs, some mathematical thing is indeed the same. That, for example, it is  interesting to extract, from behind what would be then only its masks, the computations of a certain quantity, say the ratio of a circle's circumference to its diameter,  before the quantity may have been baptised $\pi$, before it has been described as being a number or  a well-defined ratio, even before any relation has been established between it and the computation of the area of a disk.  In more sophisticated versions, by telling the story of a series of past events which have led to finally define an object in our present, history of mathematics may convey the idea that the series of past objects associated with these events were already more or less the same in the past, and have only been subsumed under the same name. That there is an identity to be detected through or behind contingent disguises concerns not only numbers or simple statements but also whole domains like algebra or statistics, or advanced concepts like group or methods of proof. In such a history, the present of mathematics speaks for its past on various scales: besides the present formulation of a statement or an object, it also defines the implicit norms for mathematical activities, for instance that they should involve proofs (even if unwritten), or that disciplines are more or less fixed. Debates of course may have been launched and corrections been made about the point of departure, the real origin: did Euclid's \emph{Elements} or Babylonian tablets contain algebra, should we begin the history of algebra only with al-Khw\=arizm\=\i\ or Fran\c cois Vi\`ete? But there were rarely doubts raised about the relevance of the question itself and the linear character of the history constructed under the key premise of an identifiable thing running all the way long. Many would agree with Andr\'e Weil that: ``More often than not, what makes [history] interesting is precisely the early occurrence of concepts and methods destined to emerge only later into the conscious mind of mathematicians; the historian's task is to disengage them and trace their influence or lack of influence on subsequent developments'' (\cite[231-232]{Weil1978}).

There are several good reasons to adopt such a point of view. For one, it is close to the ``spontaneous history'' of the working mathematician, the chronology often given at the beginning of a mathematical article in order to motivate the results it establishes.\footnote{The expression ``spontaneous history'' along with Louis Althusser's ``spontaneous philosophy'' is now often used in history and philosophy of science to designate casual remarks made by working scientists about the past of their work (resp. their general views on mathematics), see for instance \cite{Rheinberger1994}.} The topics of such a history are also more easily those of primary interest to mathematicians (recent concepts or theorems, in particular). That mathematics deals with long-lived objects may also help to consolidate its status, at times when the importance of mathematics is questioned or the population of students in mathematics is decreasing; that ``mathematical truths have been called eternal truths \dots [because] in very different expressions, one can recognise the same truths,'' in Hieronymus Zeuthen's terms (\cite[17]{LutPur1994}), guarantees a particular value for the discipline as a whole. Such a conviction may also reinforce (or be reinforced by) Platonist views of mathematics, and as such has been taken over by some philosophers. For instance, the philosopher of mathematics Jacques-Paul Dubucs, after pointing out the differences between two presentations of a proof that there are infinitely many prime numbers, one in Euclid, one modern by Godfrey H. Hardy and Edward M. Wright, claimed that proper epistemological investigation should focus on what is perceived as a stable and identical core and in particular has ``no reason to discriminate the two texts which propose \emph{fundamentally} the same proof,'' emphasising his agreement on this issue with the authors of the modern text (\cite[41]{Dubucs1991}).

\section{Ephemeral configurations}
However, William Aspray and Philip Kitcher noticed in 1988 that ``a new and more sophisticated historiography has arisen [\dots] This historiography measures events of the past against the standards of their time, not against the mathematical practices of today'' (\cite[24-25]{AspKit1988}). Indeed, innovative approaches in the history of mathematics of the last decades have often expressed misgivings over a cleaned-up history, based on a too-rapid identification of a concept or a problem, and over its historiographical consequences. In Thomas Hawkins's words: ``The challenge to the historian is to depict the origins of a mathematical theory so as to capture the diverse ways in which the creation of that theory was a vital part of the mathematics and mathematical perceptions of the era which produced it'' (\cite{Hawkins1986}). Consequently, the focus of recent history of mathematics has been much more on localised issues, short-term interests and ephemeral situations, on ``the era which produced'' the mathematics in question; and moreover it has centred on diversity, differences and changes.

Confluent factors are here at stake. One has been largely advertised. It is linked to contemporary debates in the history of science in the large and comes with the wish to take into account social aspects of mathematics and ``how they shape the form and the content of mathematical ideas'' (\cite[25]{AspKit1988}), while dimming the line between so-called internal history (that of concepts and results) and external history (that of institutions or scientific politics). Given the quantity of recent historical writing on these issues, I shall only mention a few examples. The unification of Italian states during the nineteenth century and the cultural \emph{Risorgimento} which accompanied it favoured a flourishing of mathematics, in particular a strong renewal of interest in geometry in all its forms, with Luigi Cremona, Corrado Segre, Guido Castelnuovo or Eugenio Beltrami and their followers (\cite{Bottazzini1994, BotNas2013, CCGGMV2016}). The Meiji Restoration in Japan witnessed a multifaceted confrontation between the then extremely active, traditional Japanese mathematics (\emph{wasan}) and its Western counterparts  (\cite{Horiuchi1996}). The First World War, a ``war of guns and mathematics'' as one soldier described it, did not just kill hundreds of mathematicians on the battlefields, among many millions of others: it also launched entire domains on a vast new scale, such as fluid mechanics or probability theory, and completely reconfigured international mathematical exchanges (for instance fostering a development of set theory, logic and real analysis in newly independent Poland) (\cite{AubGol2014}). One might also think of the variety of national circumstances which preceded the creation of mathematical societies in the late nineteenth and early twentieth centuries (\cite{Parshall1995}) or the various reforms in mathematical education (\cite{GisSch2011, KarSch2014}). 

At a smaller scale, specific opportunities at specific times, putting mathematicians in close contact with certain milieux, have hosted particular, sometimes unexpected, mathematical work, be it analysis in administrative reforms (\cite{Brian1994}), number theory in the textile industry (\cite{Decaillot2002}) and in computer hardware (\cite{MoBu2008}), or convexity in the military (\cite{Hoff2002}). In such cases (and in many others studied in detail over the last decades), it is not a question of superficial analogies or obvious applications; very often the ways these connections were made, the concrete manner of transmission of knowledge through personal or institutional links, the objectives pursued, are what provides impulse to a mathematical investigation, explains the formulation it takes or the particular balance between computations and theory it displays (\cite{Hoff2010, Tournes2012}). This is, by definition, ephemeral in the sense that it implies links with social situations which have their own time scale and are most certainly not ``eternal truths.''

However, another reason for this shorter-term focus has recently become even more decisive, to wit, a more acute sensibility to the multi-layered structure of mathematics and the need to study more carefully its variegated components. At first sight, it seems simple enough: mathematics is most often inscribed in texts---although ethnomathematicians also study it directly in strings or sand (\cite{Ascher1991, Ambrosio2000, PeVa2014})---, it uses words, symbols and drawings, it defines or studies certain objects, it states results and justifies them. It could seem that all we have to do is to decipher the texts and explain the objects, the results and the proofs. Of course, that there were at times debates among mathematicians about certain proofs or objects is well-known: even as late as early modern times, some would not accept a proof based on an algebraically-expressed relation and required geometric proofs in the Euclidean style, as more solid; others dismissed proofs by contradiction or, later, non-constructive proofs; the legitimacy of negative numbers or of functions without derivatives or of sets has been put into question. And around 1900, how to found mathematics---on axiomatics, on integers, or on logic---was a topic of heated controversy among mathematicians, which in turn has been studied by historians.

But more recently, other aspects have been explored, aspects which are not necessarily linked to public and noisy debates, but are part and parcel of ordinary mathematical activities. What is or could be an object or a result, how are they chosen or defined or justified, has changed in time; it has also been seen as depending on the place or the author. What sort of question is considered interesting, by whom and why; which criteria are required to make an argument convincing or a solution satisfactory, again for whom and why; all these aspects and their relations are worth being studied for their own sake. Historians have, for instance, shown that arguments in words or symbols may rely upon, or be inspired by, or sometimes even been replaced by figures, diagrams, tables, instruments.\footnote{For illustrations of these different cases, see for instance \cite{Netz1999, Lorenat2014, MPS2012, CCFR2003, Durand2010, TobTou2011, FlaNab}. References being too numerous for exhaustivity, those given are only illustrative.} That an acceptable answer to a problem may be, at times, and for certain groups of mathematicians, a single number, an explicit description of all the solutions, an equation, an existence theorem, or the creation of a new concept (\cite{Goldstein2001, Chorlay2010, Ehrhardt2012}). 
The way various domains are defined and interact, or are perceived as distinct, has also changed within mathematics, but also between mathematics and other domains, in particular physics (\cite{Archibald1989, Gray1999, SchSch2011}).

What are called ``epistemic values,'' that is the internalised criteria of what constitutes good mathematics at one time, have also been studied: rigour is the most obvious perhaps and has a complex history, but universalism or effectiveness or generality or naturalness may appear at some moments to be even more decisive (\cite{Mehrtens1990, Rowe1992, Schubring2005, Corry2004, Gene2016}). 
Working mathematicians have usually their own answer to these questions, but the point here is to reconstruct the whole range of positions at a given time, in a given milieu, and to understand their effect in mathematical work. For instance, in the dispute between Leopold Kronecker and Camille Jordan in 1874 about what we now see as the \emph{same} reduction theorem for matrices, differences in formulation (elementary divisors on one side, canonical forms on the other), and in disciplines (invariant theory vs group theory) were at play, but also differences in conceptions of generality (\cite{Brechenmacher2016}). 

Last, but not least, the way mathematics is made public and circulates has been proved to be both a serious constraint on its form and an essential factor in its transmission; the organization of correspondence among mathematicians (indeed the mere form of a mathematical letter), the creation of mathematical research journals in the nineteenth century and their different organisation through time, the advent of academies, seminars and conferences, teaching programs and textbooks, for instance, have all been scrutinised (\cite{Peiffer1998,Rowe2004,Schubring1985a,Schubring1985b,AusHor1993,Verdier2009,Gerini2002,RemSch2010}). The last two aspects could also be considered as links between mathematical developments and general cultural issues, but here the emphasis is on the combination of these components inside mathematical texts themselves. 

Depending on one's own tastes, this sheer variety in the course of time may appear fascinating or an irrelevant antiquarian interest. However, we now have enough evidence that all these aspects may count for understanding the development of mathematics. Mathematics weaves together objects, techniques, signs of various kinds, justifications, professional lifestyles, epistemic ideas. Recent biographies, indeed, offer successful examples of the study of such articulations (\cite{Parshall2006, Crilly2006, Alfonsi2011}). But historians have also studied these components separately,  in a comparative way, in order to display their range and their evolution in time. These components have distinct time-scales and changes do not occur simultaneously. Even when one is able to understand a long-term development of one component (for instance, of mathematical publishing), its articulation with other components is generally stable only over a shorter period. 
A further difficulty is that if concepts or theorems are aspects of which the mathematician is aware (and very much so), some of the components I mentioned are much more implicit, or are operational at a collective, not at an individual, level; they can be best detected and analysed for an entire group (\cite{Goldstein1999}). All this explains the interest in studying what I am calling here ``configurations,'' a word borrowed from the sociologist Norbert Elias. Elias wanted to set himself apart from previous sociological theories based on an a priori hierarchic opposition between individuals and society, and he promoted the idea of first studying configurations formed by interactions between persons in interdependence, at different scales, be they players in a game or workers in an enterprise. For us, configurations organise texts and persons, coordinating some of the components we have mentioned (we shall see concrete examples later).

A last reason has favoured more localised studies: a critical outlook by historians of mathematics on their own practice. Words like ``discipline'' or ``school'' have been used in the past without further ado, in particular because they were terms inherited from mathematicians themselves. Recent work has shown that to define and use them more carefully gave a better grasp for describing the past. For instance, using a characterisation of a discipline as a list of internal elements (core concepts, proof system, etc), Norbert Schappacher and I were able to distinguish in the lineage of C. F. Gauss's \emph{Disquisitiones Arithmeticae} those parts which fused into a existing discipline (his treatment of the cyclotomic equation, for instance, which had a potent effect on the theory of equations), those parts which emancipated themselves as autonomous disciplines  with their own programmes and priorities (quadratic forms in the middle of the nineteenth century, reciprocity laws in the theory of number fields later), and those parts which, in the nineteenth century, were taken over in an isolated manner by some mathematicians (primality tests).\footnote{\cite{GolSch2007}, see also \cite{Gauthier2009, Gauthier2010}. On the issue of ``school,'' see the synthesis (\cite{Rowe2003}).} Caution also applies to common descriptors of historical phenomena themselves, such as ``context'' (\cite{Ritter2004}), ``\emph{longue dur\'ee}'' (\cite{AubDah2002}) or ``revolution'' (\cite{Gillies1992}). This reflexivity has also permitted historians to find counter-examples to overly-crude hypotheses on the long-term development of mathematics (\cite{GolRit2003, GilGui2015}).

Thus, in the last decades, historians have used these analyses of   configurations to deconstruct the identifications provided by mathematical works in the course of time. Examples include statements such as the fundamental theorem of algebra (\cite{Gilain1991}), Sturm's theorem (\cite{Sinaceur1991}), Fermat's theorem that  ``the area of a right-angled triangle in integers cannot be a square''  (\cite{Goldstein1995}) or the decomposition theorem of matrices (\cite{Brechenmacher2007}); concepts like ideals (\cite{Edwards1980, Edwards1992})  or points (\cite{Schappacher2010}); or even whole domains like Galois theory (\cite{Ehrhardt2012}).   As explained by Andr\'e Weil, for instance, Fermat's 1640 statement, and even his proof, can be identified and seen as \emph{the same} as a special case of the Mordell-Weil theorem, according to which the group of rational points of the elliptic curve defined by $y^2z=x^3-xz^2$ is $\mathbb{Z}/2\mathbb{Z} \times \mathbb{Z}/2\mathbb{Z}$.  But we can also see this identification as a historical problem: it requires first to reconstruct various configurations involving each of the statements,  and then to understand how they have come to be identified---in other words, to study also the mathematical \emph{work} that provided such retrospective identification.

Here,  differences are obvious, but some cases are  more delicate. This can be illustrated by the theorem that there are exactly twenty-seven lines on a non-singular cubic projective surface; since its statement (and proof) in 1849 by both Arthur Cayley and George Salmon, its formulation has remained remarkably stable for more than a century. But what changed is its association with other problems: as shown in (\cite{Le2015}), it is for instance in tandem with the fact that there are 9 inflection points on a cubic projective plane curve and other analogous statements that it played a decisive role for the assimilation of group-theoretical methods by geometers \emph{before} Felix Klein's Erlangen Program; this specific configuration of questions and disciplinary issues, around the so-called ``equations of geometry,'' lasted only a few years, but was a key feature in the transmission of the theorem. 

In what follows, I would like to illustrate these issues with what could be described as a minimal example:  the concrete case of a rather technical and apparently stable concept, that of an Hermitian form,\footnote{For our purpose, such a form will be simply an expression of the type $\sum_{i, j=1}^n a_{ij}x_i \overline{x_j}$, with coefficients $a_{ij}$  in $\mathbb{C}$, such that $a_{ji}=\overline{a_{ij}}$ (here, the bar designates the complex conjugation); in particular, the diagonal coefficients $a_{ii}$ are real numbers.}  over a  short period of time,  the second half of the nineteenth century.

\section{Did Hermite invent Hermitian forms?}
From all accounts, for instance by \cite[612-613]{Vahlen1900} or \cite[vol.~3, p.~269]{Dickson1919}, forms of this type first publicly appear in an article authored by Charles Hermite and it is thus legitimate to ask if Hermite indeed invented Hermitian forms, and how.

At the very beginning of this paper (\cite{Hermite1854}), Hermite explained:

\begin{quote}
One knows how easily one can extend the most fundamental arithmetical concepts coming from \emph{real integers} to complex numbers of the type $a+b\sqrt{-1}$. Thus, starting from elementary propositions concerning divisibility, one quickly reaches those deeper and more hidden properties which rely upon the consideration of quadratic forms, without changing anything essential in the principle of methods which are proper to real numbers. In certain circumstances, however, this extension seems to require new principles and one is led to follow in several different directions the analogies between the two orders of arithmetical considerations. We would like to offer an example to which we have been led while studying the representation of a number as a \emph{sum of four squares}.
\end{quote}

No references were provided in this introduction, but the allusions seem rather clear, both then and now. Carl Friedrich Gauss, searching for an extension to higher powers of the reciprocity law for squares that he had proved in his 1801 \emph{Disquisitiones Arithmeticae}, launched an arithmetical study of what he called ``complex integers'' (now Gaussian integers), that is ``complex numbers of the type $a+b\sqrt{-1}$,'' with $a$ and $b$ ordinary integers (\cite{Gauss1831}). Among them he defined prime complex numbers and units, proved the factorisation of the ``complex integers'' into a product of these prime numbers (unique up to units and to the order of the factors), showed how to extend Euclidean division and  congruences to these numbers: in short, ``the elementary propositions which concern divisibility.'' In 1842, Peter-Gustav Lejeune-Dirichlet  began to study  ``those deeper and more hidden properties which rely upon the consideration of quadratic forms,'' in particular the representation of Gauss's ``complex integers'' by what will come to be known as Dirichlet forms at the end of the century, that is, quadratic forms $f(x,y) = ax^2+2bxy+cy^2$, where the coefficients $a, b, c$, and eventually the values taken by the indeterminates $x$ and $y$, are also ``complex integers'' (\cite{Dirichlet1842}). 

Such a discontinuous chronology (1801--1831--1842--1854) might lead us to the topic of Hermite's paper, but would not be sufficient to understand Hermite's background or  point of view. First of all, Gauss's discussion of his complex integers appeared in Latin in the proceedings of the G\"ottingen Society of Science, with limited distribution. But as early as 1832, Dirichlet explained Gauss's work and completed it in what was at the time the only important journal entirely devoted to mathematics, August Leopold Crelle's \emph{Journal f\"ur die reine und angewandte Mathematik}, created in 1826. Dirichlet, who had spent several years in Paris, was an important go-between for mathematics: his 1832 article, written in French, was clearly aimed at an international audience. In the same decade, he would use new tools developed by analysts, in particular Fourier series, to complete proofs of Gauss's statements and revisit a number of his arithmetical results, stressing their links to various areas of mathematics in a way that would draw greater attention to them.  In 1840, for instance, a letter to Joseph Liouville, on the occasion of a French translation of one of his papers,  announced to the French community his current interest for ``extending to quadratic forms with complex coefficients and indeterminates, that is, of the form $t+u\sqrt{-1}$, the theorems which occur in the ordinary case of real integers. If one tries in particular to obtain the number of different quadratic forms which exist in this case for a given determinant, one arrives at this remarkable result, that the number in question depends on the division of the lemniscate; exactly as in the case of real forms with positive determinant, it is linked to the division of the circle'' (\cite{Dirichlet1840}). The lemniscate pointed to the integral $\int \frac{dx}{\sqrt{1-x^4}}$ and to elliptic functions, then at the forefront of research, and which will be soon the main topic of interest of the young Hermite. 

Then, the ten years preceding Hermite's 1854 paper were  turbulent years for complex functions and numbers, and to a lesser extent, for quadratic forms. In a footnote of his 1832 paper, Dirichlet  announced, somehow optimistically, that numbers of the type $t+u\sqrt{a}$, for $t$ and $u$ integers and $a$ an integer without square divisors, give rise to theorems analogous to those on Gaussian integers, and with similar proofs. In 1839 (with a French translation in Liouville's \emph{Journal de math\'ematiques pures et appliqu\'ees} three years later), Carl Jacob Jacobi, again recalling Gauss's theory of complex integers, showed that a prime number $p=8n+1$ can be written  as a product of \emph{four} complex numbers, each of them a linear combination with integral coefficients of powers of a given $8^{\textrm{th}}$-root of unity, such that the three decompositions of $p$ as $x^2+y^2$, $x^2+2y^2$ and $x^2-2y^2$ be issued ``from a common source.'' Announcing similar results for a prime $p=5n+1$ (and $5^{\textrm{th}}$-roots of unity), but with no hint of a proof, Jacobi provided the spur for decisive work by several younger mathematicians in the 1840s. These included Gotthold Eisenstein's approach to complex multiplication, and Ernst Eduard Kummer's theory of ideal numbers (in what we call now cyclotomic rings); Kummer's display that unique factorization failed in general certainly crushed Dirichlet's 1832 hopes and showed, as Hermite pointed out in 1854, that ``in certain circumstances, this extension [of arithmetic to complex numbers] seems to require new principles.'' It also included Hermite's own first research on quadratic forms, concerning which he wrote directly to Jacobi from 1847 on.\footnote{This is explained in more detail in \cite[39-51]{GolSch2007}.} Personal relations here reinforced the circulation of the articles; Hermite was informed of Kummer's approach to the arithmetic of complex numbers by the mathematician Carl Wilhelm Borchardt  during the latter's 1847 Parisian tour and he met Dirichlet and Eisenstein, among others, during his own trip to Berlin at the beginning of the 1850s.

Hermite's work on forms arose at least as much from his close reading of Jacobi---that on the decomposition of primes, but also that on elliptic functions---as from Dirichlet's articles on complex numbers. In his letters to Jacobi on quadratic forms, \cite{Hermite1850}  considered forms with any number of indeterminates and real coefficients (instead of the two indeterminates and integral coefficients of Gauss's \emph{Disquisitiones}). His main theorem was to establish that there exists a (non-zero) value of the form, when evaluated on integers, which is less than a certain bound, depending only on the number of indeterminates and on the determinant of the form, but not of its coefficients.\footnote{It was common at the time not to distinguish explicitly among indeterminates, variables and values; or to give general statements without a clear list of exceptions. Moreover  different conventions coexisted, for instance for the definition and the sign of  determinants, etc. For reasons of space, and although these questions may be revealing and have been taken into account in historical work, I shall not in general discuss them here.} That is:

Let $f(x_0,x_1\ldots, x_n)$ a definite positive quadratic form with $n+1$ indeterminates, \emph{real} coefficients, and determinant $D$. Then, there exist $n+1$ \emph{integers} $\alpha$,
$\beta$, $\ldots$, $\lambda$, such that

\begin{equation} \label{main}
0 < f(\alpha, \beta, \ldots, \lambda)< (\frac{4}{3})^{n/2} \root {n+1}\of {\mid D\mid}.
\end{equation}

Although the formulation does not make it obvious, this statement is closely related to the classification of forms. For definite binary forms, $f(x,y)=ax^2+2bxy+cy^2$, with determinant $D=ac-b^2$ and integral coefficients,  Gauss's reduction theory stated that among all the forms arithmetically equivalent to a given form $f$ (that is all the forms which are derived from $f$ by an invertible linear transformation of the indeterminates  with integral coefficients), there exists one, say $F$, whose first coefficient $A$ is less than $2\sqrt{\frac{\mid D\mid}{3}}$. This coefficient $A=F(1,0)$ is a value of $F$, and thus it is also (by linear transformation) a value at  some integers of all the equivalent forms $f$. For these forms,  this is exactly Hermite's main theorem for $n=1$. 
Hermite used his general statement to prove that if one restricts oneself to forms with integral coefficients, the number of classes of arithmetically equivalent forms for a given determinant is finite, and in particular that, for $n=1, 2, \cdots, 6$, there is only one class of forms with $\mid D\mid=1$, represented in each case by the sum of $n$ squares. 

Classification was a key objective for Hermite: ``[T]he task for number theory and integral calculus,'' he wrote to Jacobi, ``[is] to penetrate into the nature of such a multitude of entities of reason, to classify them into mutually irreducible groups, to constitute them all individually through characteristic and elementary definitions'' (\cite[286]{Hermite1850}). The striking echo of a quasi-botanical project in this quote is to be taken seriously: for Hermite 
 and some of his contemporaries, the emphasis on classification directly came from a  view of mathematics as a natural science (\cite{LP2016}). ``Collecting and classifying'' was also a very strong incentive for invariant theorists like Arthur Cayley (\cite[193-195]{Crilly2006}), and it was not limited to them, nor to the 1850s: in 1876, still, Leo K\"onigsberger wrote for instance: ``It seems to me that the main task now just as for descriptive natural history consists in gathering as much material as possible and in discovering principle by classifying and describing this material'' (File H1850(6), Staatsbibliothek zu Berlin, Handschriftenabteilung).

Hermite's main theorem was quite versatile: with it, for instance, Hermite simultaneously showed how to approximate real numbers by rationals and  proved two statements left unproved by Jacobi, that there is no complex function of one variable with three independent periods and, as announced above, that prime numbers of the form $p=5n+1$ can be decomposed into a product of 4 linear combinations with integral coefficients of powers of a fifth root of 1. Characteristically, the unity of mathematics is here found in  the bridges between analysis, algebra and arithmetic. In each case, the whole point is to choose correctly a quadratic form (or sometimes a continuous family of them) to encapsulate the phenomenon under scrutiny and to combine inequalities provided by the theorem with integrality properties.

Hermite also adapted this construction to discuss the divisors of forms of the type $x^2+Ay^2$, and then, in 1854, took the natural step of testing it on the famous theorem that every integer is the sum of four squares. 

Let $A$ be a non-zero integer (without loss of generality, one may assume that 4 does not divide $A$). As a preliminary step, Hermite first showed how to find integers $\alpha$ and $\beta$ such that $\alpha^2+\beta^2 \equiv -1 \bmod (A)$.

He then introduced the quadratic form with 4 variables:
\begin{equation*}
f(x,y,z,u)=(Ax+\alpha z+\beta u)^2+(Ay-\beta z+\alpha u)^2+z^2+u^2,
\end{equation*}
\noindent with determinant $A^4$. When $x$, $y$, $z$, $u$ are integers, the value $f(x,y,z,u)$ is also an integer and, thanks to the choice of $\alpha$ and $\beta$, it is divisible by $A$.  On the other hand, Hermite's main theorem states that there exist integers $x, y, z, u$ such that $0<f(x,y,z,t)<(\frac{4}{3})^{\frac{3}{2}}\sqrt[4]{A^4}$, that is, $0<f(x,y,z,u)< 1.54  A$. The only possibility is that $f(x,y,z,u)=A$, which expresses $A$ as a sum of 4 squares.

However, with classification in mind, Hermite decided to reformulate this proof: instead of $f$, he considered $\frac{1}{A} f$ (the determinant of which is then 1). It can be written, when rearranging terms, as

\begin{multline}\label{quatre}
\frac{1}{A}f(x,y,z,u)=A(x^2+y^2) + 2\alpha (zx+yu)+2 \beta (xu-zy)\\+\frac{\alpha^2+\beta^2+1}{A}(z^2+u^2).
\end{multline}

As said above, Hermite had already proved that, up to arithmetic equivalence,  there is only one quaternary quadratic form of determinant 1 with integral coefficients, the form $X^2+Y^2+Z^2+U^2$. That is, there exist integers $m, m', \cdots, n, n', \cdots $ such that the change of variables
\begin{align*}
X &= mx+m'y+m''z+m'''u \\
Y &= nx+n'y+n''z+n'''u \\
Z &= px+p'y+p''z+p'''u \\
U &= qx+q'y+q''z+q'''u
\end{align*}
\noindent transforms $X^2+Y^2+Z^2+U^2$ into the form $\frac{1}{A}f$. By identification of the term in $x^2$, for instance, one obtains $A=m^2+n^2+p^2+q^2$, as desired. 

Viewing mathematics as a natural science also entailed the highlighting of specific practices.  ``[My own work] would very strikingly illustrate how much observation, divination, induction, experimental trial, and verification, causation, too [\dots] have to do with the work of the mathematician,'' claimed James Joseph Sylvester at the British Association for the Advancement of Science in 1869 (\cite{Parshall2006}) and Hermite himself repeated several times that ``the most abstract analysis is for the most part an observational science'' (\cite[403]{Goldstein2007}). Observed carefully, the very shape of the form $\frac{1}{A}f(x,y,z,u)$ (equation \ref{quatre}) suggests ``complex numbers.'' At this stage, indeed, Hermite introduced a new set of indeterminates: 
\begin{align*}
v=x+\sqrt{-1}y &\qquad V=X+\sqrt{-1}Y\\
v_0=x-\sqrt{-1}y &\qquad V_0=X-\sqrt{-1}Y\\
w=z+\sqrt{-1}u &\qquad W=Z+\sqrt{-1}U\\
w_0=z-\sqrt{-1}u &\qquad W_0=Z-\sqrt{-1}U
\end{align*}

\noindent and \emph{restricted} the type of linear transformations to be taken into account.

\begin{quote}
One can distinguish among all real transformations between the two groups of variables $x, y, z, u$ on one hand, and $X, Y, Z, U$ on the other, those which can be expressed as:
\begin{align*}
v&=aV+b W\\
v_0&=a_0V_0+b_0 W_0\\
w&=c V+d W\\
v_0&=c_0V_0+d_0 W_0\\
\end{align*}
\noindent where $a, b, c, d$ are arbitrary imaginary numbers and $a_0, b_0, c_0, d_0$ their conjugates. Thus one obtains a perfectly defined class of real transformations.
\end{quote}

In a modernised matrix notation, these transformations can be expressed as
special transformations with 4 real variables:
\[\begin{pmatrix}
x\\
y\\
z\\
u
\end{pmatrix}
=
\begin{pmatrix}
   \Re{(a)}  & -\Im{(a)}  &  \Re{(b)}  & -\Im{(b)} \\
     \Im{(a)}  & \Re{(a)}  &  \Im{(b)}  & \Re{(b)}\\
       \Re{(c)}  & -\Im{(c)}  &  \Re{(d)}  & -\Im{(d)}\\
        \Im{(c)}  & \Re{(c)}  &  \Im{(d)}  & \Re{(d)}\\
\end{pmatrix}
\begin{pmatrix}
X\\
Y\\
Z\\
U
\end{pmatrix}.\]

These transformations, in turn, preserve specific \emph{quadratic} forms that Hermite introduced at this point (those we now call Hermitian forms): 
\begin{equation}\label{herm}
f (v,w)=Avv_0+Bvw_0+B_0wv_0+Cww_0
\end{equation}

\noindent where $A$ et $C$ are real numbers, and $B$ a complex number ($B_0$ its complex conjugate). 

If one replaces the complex variables $v$ and $w$ by their real and imaginary parts, that is by real variables, one finds: 
\begin{equation*}
f(x,y,z,u)=A(x^2+y^2) + 2\Re(B) (zx+yu)+2 \Im(B) (xu-zy)+C(z^2+u^2),
\end{equation*}

\noindent of which the form (\ref{quatre}) used to prove the theorem of four squares is indeed a prototype. In Hermite's terms:

\begin{quote}
Considered with respect to the  original variables $x, y, z, u$, these forms are entirely real, but their study, with respect to the transformations we have defined previously, essentially relies upon the use of complex numbers. One is then led to attribute to them a mode of existence singularly analogous to that of binary quadratic forms, although they essentially contain four indeterminates (\cite[346]{Hermite1854}).
\end{quote}

Moreover, if the linear transformations previously considered have a determinant of complex norm 1, they leave invariant the quantity $\Delta=BB_0-AC$, which plays the role of the determinant for the form (\ref{herm}). Hermite then undertook a classification of these forms.

Two important features should be  underlined here. The first one is the central role played by the linear transformations, as key tools for classifications: this is the restriction on these transformations that defined the new type of forms. In several memoirs around 1850, Hermite and some of his contemporaries,  Cayley, Borchardt, Eisenstein, Otto Hesse, for instance, described several types of equivalence among forms, depending on the transformations which were taken into account: either they looked for transformations keeping a given form or function invariant, or, fixing a group of transformations, they studied the forms which are left invariant by them.  The classification of forms with real coefficients (through linear transformations with real coefficients), which led to the notion of signature, belonged to the  same programme. In 1855, in a study on the transformations of Abelian functions with 2 variables (more precisely of their periods), Hermite again introduced other ``particular forms with 4 indeterminates, where one does not use as analytical tool the most general transformations among 4 variables, but particular transformations \dots which reproduce analogous forms'' (\cite[785]{Hermite1855a}). These forms and transformations would be one of the origins of the symplectic forms and group (\cite{Brouzet2004}).

This point of view, changing the group of transformations operating on the variables in order to delineate which type of forms or functions will be studied, was at the time tightly linked to invariant theory (\cite{Parshall1989}) and its applications were varied. It is in the context of Sturm's theorem on the number of roots of an algebraic equation that belong to a given domain that Hermite generalised his 1854 construction to quadratic forms with $2n$ ``pairwise conjugate'' indeterminates (again the index 0 designates the complex conjugation), 

\begin{equation*}
f(x_1, x_2, \dots, x_n, x_{1,0}, \dots, x_{n,0}) = \sum_{i, j} a_{i, j} x_i x_{j,0}, 
\end{equation*}

\noindent with $a_{i,j}$ and $a_{j, i}$ complex conjugates (thus $a_{i,i}$ real numbers) (\cite{Hermite1855b, Hermite1856}); he called them simply ``quadratic forms with conjugate imaginary indeterminates.''

The second  feature worth stressing concerns the way mathematical objects are introduced. Hermite's ``forms with conjugate imaginary indeterminates''  were for him quadratic forms of a specific type, not a new type of objects defined in an ad hoc way, or by simple analogy, to accommodate complex numbers. They are  distinguished \emph{among} a larger family of well-known objects because of their special properties (here their stability under a certain group of transformations), as a new species could have been. In Hermite's view, coherent with that of mathematics as a natural science,  mathematicians do not, should not, create their objects: they ``meet them or discover them and study them, like physicists, chemists and zoologists'' (\cite[vol. 2, p. 398]{HerSti1905}). Hermite's key role in a history of Hermitian forms cannot be doubted, but its description is thus delicate. Besides the problems raised by the word ``invention'' in this particular context, another issue  is directly connected to our main point, that of sameness: Alfred Clebsch in 1860 and 1863 (\cite{Clebsch1860, Clebsch1863}),  Elwin Bruno Christoffel a bit later (\cite{Christoffel1864}) also introduced, in completely different contexts, mathematical objects close to Hermite's one: Clebsch studied square \emph{matrices} with complex entries, such that entries symmetric with respect to the diagonal are complex conjugates (that is, now, Hermitian matrices), Christoffel considered \emph{bilinear forms} on two sets of indeterminates, $\phi =\sum [gh]u_gv_h$ and their values when $u_g$ and $v_g$ are set to be complex conjugate numbers, under the assumption that the coefficients $[gh]$ and $[hg]$ are complex conjugates.  Both referred to Hermite's 1855 article (only as an afterthought, in 1864, for Clebsch) for specific results, but none  identified his  own construction to Hermite's one nor gave any special name to it. These various points of view on Hermitian forms will be unified much later in the century. 

The configuration about the appearance of Hermitian forms I have  briefly sketched thus includes local incentives, a series of specific mathematical themes, a collection of objects and the disciplinary tools available to study them (here for instance the reduction theory of forms), the state of the art on certain topics (for instance complex numbers), but also here a new emphasis on linear transformations and a mathematical world-view with an impact on the practice of mathematics and how an entity is constructed and accepted. None of these components is proper to Hermite alone, even if their coordination may be, and they would require to be studied as collective phenomena.  Some of these components  will evolve separately and at different rates, losing their  connections with the development of our forms ``with conjugate imaginary indeterminates.''  

Tracing, on the other hand, the fate of this particular type of forms is not obvious. In 1866, for instance, Cayley discussed ``Hermitian matrices'' and associated forms, but they were those attached to the transformation of Abelian functions, not ``our'' forms and transformations (\cite{Cayley1866}).  The variety of names or of notations (for complex conjugation, in particular) also forbids us to use simple criteria or visual display. 

A technique which has recently demonstrated its effectiveness for other topics (\cite{Goldstein1999, Brechenmacher2007, Roque2015}) is first to systematise the search for relevant writings in review journals (such as the \emph{Jahrbuch \"uber die Fort\-schrit\-te der Mathematik}), then to reconstruct their own networks of references in order to track down how Hermitian forms circulated after 1854. Thamous, an electronic tool  developed by Alain Herreman for collaboratively constructing relational databases, was used here extensively to  register cross-references and other links between articles and to organise the corpus. For reasons of space, and with a view to the long-term question, I shall only briefly report on some of the configurations in which our forms occur before the First World War. 

\section{Picard's and Bianchi's groups}

Both \'Emile Picard and Luigi Bianchi were born in 1856, at more or less the same time as Hermitian forms. And both would be instrumental in bringing them back center-stage, this time with a group-theoretical and geometrical apparatus. Although some of their results are very close (even leading to some tensions), their backgrounds are quite different. Picard is best known for his work on complex analysis and its application to algebraic surfaces (\cite[245]{Houzel1991}) and, more generally, for extending to dimension 2 a number of results first established in dimension 1. At the beginning of the 1880s, for instance, he proved, in parallel with the 1-dimensional case, that of the elliptic curve, that surfaces which can be parametrised by Abelian functions have (under some restrictions on their singularities) a geometric genus $p_g\leq1$ (\cite{Picard1881}). 

But at the very same time and in close proximity, a great enterprise was underway; Henri Poincar\'e had just developed his theory of ``Fuchsian'' and ``Kleinian'' functions (now both considered particular cases of automorphic functions), that is, meromorphic functions of one complex variable $z$ inside a certain disk such that 

\begin{equation} 
f(\frac{az+b}{cz+d}) =f(z), 
\end{equation}

\noindent where the Moebius transformations $z \rightarrow \frac{az+b}{cz+d}$ belong to a  discontinuous group of invertible transformations; the groups and the functions were ``Fuchsian'' for Poincar\'e when the coefficients $a, b, c, d$ are real.  Poincar\'e had interpreted the Moebius transformations on the unit disk (or equivalently on the upper half plane) as isometries in a non-Euclidean, hyperbolic, geometry and he had constructed fundamental domains for the Fuchsian groups. He had also shown, given a Fuchsian group, how to construct Fuchsian functions invariant under the group and how to connect them to the solution of differential equations (\cite{Gray1989}). 

Picard was at first in search of a first, analogous, example with two variables. For this he introduced the family of curves of equation $z^3=t(t-1)(t-x)(t-y)$, with $x$ and $y$ two real parameters, and he studied the periods of $\int z^{-1} dt$ (\cite{Picard1882a}). These periods, as functions of $x$ and $y$, satisfy a system of partial differential equations, which admits a basis of  three independent solutions $A_1, A_2, A_3$. The functions $u=\frac{A_2}{A_1}$ and $v=\frac{A_3}{A_1}$ of  the two variables $x$ and $y$ can be inverted, providing Picard with what he was looking for, two uniform functions of the two variables $u, v$, defined on the domain $2 \Re(v)+ \Re{u}^2+\Im(u)^2 <0$. He also computed a group of linear transformations under which his functions would remain invariant. But, as he explained later, in order to obtain more than this one isolated example, ``the thought of a recourse to ternary quadratic forms with conjugate indeterminates showed [him] the way out''  (\cite[36]{Picard1889}).

A ternary form with conjugate indeterminates, analogous to (\ref{herm}), 
\begin{multline}\label{pic}
f (x, y, z, x_0, y_0, z_0)=axx_0+a'yy_0+ a''zz_0+ byz_0+b_0y_0z\\+ b'zx_0+b'_0z_0x+ b''xy_0+b''_0x_0y
\end{multline}
\noindent with real coefficients $a, a', a''$, and the index $0$ designating the complex conjugation, as for Hermite, becomes becomes either $\pm (xx_0+yy_0+zz_0)$ or $\pm (xx_0+yy_0-zz_0)$ by an adequate linear transformation of the indeterminates. The first case corresponds to definite forms, the second to indefinite forms. In the 4-dimensional (real) space defined by the two complex variables $u=\frac{x}{z}, v=\frac{y}{z}$, the equation $f=0$ represents  a 3-dimensional (real) hypersurface (for instance, the form $xx_0+yy_0-zz_0$ corresponds to the hypersphere $uu_0+vv_0 =1$). Picard then restricted the coefficients of the form to be ``complex integers'' \footnote{It means, for most of his work,  Gaussian integers. Picard did in fact consider complex quadratic numbers in general, but his definition of ``integers'' here was not correct.} and studied the linear transformations  leaving the form invariant or, more precisely, their non-homogeneous versions operating on $u$ and $v$.    If the form is definite, the group that one obtains is finite,  but Picard's hope  was to use the groups obtained in the indefinite case as  analogues of  Fuchsian groups. He had first to justify their existence and main properties (\cite{Picard1882b}).
 
This was done by showing how the interior of the hypersurface $f=0$ can be cut into a tessellation of polyhedra, transformed  one from another by an infinite discontinuous group. The procedure corresponded geometrically to Hermite's so-called continuous reduction process for quadratic forms:  Picard associated to the indefinite form a 2-parameter family of definite forms, for which a well-known theory of reduction existed and could be expressed as in the classical binary case by conditions of inequalities on the coefficients of the forms. To each choice of the parameters, that is, to each definite form in the family, Picard could associate  in a 4-dimensional space a point, which belongs to a specific domain; this domain is a fundamental polyhedron defined by the inequalities of the reduction theory, exactly if the form is reduced. By changing the parameters continuously, the associated definite form may cease to be reduced, but, as in the binary classical case, can be transformed into an equivalent reduced form by a linear transformation, and correspondingly the point associated to it is transformed into a point in the fundamental polyhedron. The initial ternary indefinite form thus gave rise to an infinite discontinuous group of transformations (that Picard called hyperfuchsian), associated to the tesselation of the hypersurface. In later papers, Picard constructed \emph{hyperfuchsian functions} defined in the interior of the hypersurface and invariant by such a group and showed that hyperfuchsian functions corresponding to the same group could be expressed as rational functions of three of  them, linked by an algebraic relation. Picard also developed an analogous arithmetical study of \emph{binary} forms with conjugate indeterminates, interpreting their reduction geometrically in terms of domains on the plane limited by arcs of circles and, this way, constructed afresh Fuchsian groups (\cite{Picard1883, Picard1884, Picard1891}). 

When he picked up the topic in 1890, Luigi Bianchi was coming from a quite different background, which included number theory, group theory and some geometry, but not the  analytical connections dear to Picard. Bianchi had studied with Felix Klein in G\"ottingen during his European post-doctoral tour and the main reference in his papers on the arithmetic of forms is Klein's \emph{Vorlesungen \"uber die Theorie der elliptischen Modulfunctionen}, completed by Robert Fricke, which had just appeared in 1890. As \cite[313]{Bianchi1891} explained it:

\begin{quote}
The geometrical method, on which Professor Klein bases the arithmetical theory of the ordinary binary quadratic forms, may be applied with the same success on a larger scale. To prove this is the aim of the following development which will treat in the same way the theory of Dirichlet forms with integral complex coefficients and indeterminates and of Hermitian forms with integral complex coefficients and conjugate indeterminates.
\end{quote}

Indeed, Bianchi had just studied the arithmetic of Dirichlet forms with Gaussian-integer coefficients,  in order  to complete Dirichlet's results on the number of classes of such forms. He then proceeded to complete some points in Picard's study of the \emph{arithmetic} of Hermitian forms, launching both an extension to forms with coefficients in any quadratic field and a detailed examination of the associated groups of transformations and their subgroups of finite index (\cite{Brigaglia2007}). In particular, Bianchi would handle forms whose coefficients are integers in quadratic fields $\mathbb{Q}{\sqrt{-D}}$ for $D=1, 2, 3, 5, 6, 7, 10, 11, 13, 15, 19$, displaying a good knowledge of Richard Dedekind's theory of ideals. He also extended the main group of transformations to include those whose determinant is any unity and, explicitly following an idea of Fricke, those of the type $z \rightarrow \frac{az_0+b}{cz_0+d}$ (the index, as before, indicating the complex conjugation) and he computed the corresponding fundamental polyhedra. As announced, all along, and in particular in his synthesis published in \emph{Mathematische Annalen} in German, Bianchi handled side by side the arithmetic of Dirichlet and of Hermitian forms and their geometrical interpretations, in particular the determination of their fundamental polyhedra. The emphasis moved from a configuration where complex analysis plays a key role to one centered on number theory and group theory; Hermitian forms are then number-theoretical objects parallel to those introduced by Dirichlet.

Picard's and Bianchi's results were integrated by Klein and Fricke  in their 1897 \emph{Vorlesungen \"uber die Theorie der automorphen Funktionen}, but through Bianchi's point of view. They would then be taken up and developed through different methods by a number of mathematicians in the following decades, Onorato Nicoletti, Otto Bohler, Leonard Dickson, Georges Humbert, Gaston Julia, Hel Braun, Hua Luogeng, and many others. 

\section{A theorem and three authors}
The question of the \emph{finite} subgroups of the linear groups returned to the forefront as the theorem:
``For any finite group of $n$-ary linear homogeneous transformations, there exists an $n$-ary positive definite Hermitian form, say, $\sum a_{ik} x_i\bar{x_k}$, that the group leaves absolutely invariant.''\footnote{Note that for the authors we are discussing,  a linear transformation of the group  operates on the indeterminates $x_i$, the transformation obtained by conjugation of the coefficients operating on the $\bar{x_i}$. Again, the terminology and the viewpoint vary slightly according to the authors: the word ``substitution'' instead of ``transformation'' is still widely used (see the titles of the papers), ``homogeneous'' to indicate elements of $GL_n(\mathbb{C})$ is sometimes omitted, etc.}

For once, the date of this theorem is rather precise: July 1896! The  day, on the other hand, and the author, are another story, which illustrates  the new role of mathematical societies and seminars and their interaction with journals for the circulation of mathematics, as well as the wider internationalization of mathematics at the end of the nineteenth century (\cite{Parshall1995}).
On Monday July 20, 1896, the statement appeared in a note by Alfred Loewy  presented to the French Academy of Sciences by Picard for insertion in the \emph{Comptes rendus} (\cite{Loewy1896}). Following the rule of this journal to accept only very short communications, it consisted mainly of an announcement and contained no proofs. Two years earlier, Loewy had obtained his thesis under Ferdinand Lindemann, with a work on the transformation of a quadratic form into itself. In his 1896 note, Loewy first gave a condition for a linear transformation to fix an arbitrary  bilinear form with conjugate indeterminates, $\sum a_{ik} x_ix_k^0$, with arbitrary coefficients (the complex conjugation is  indicated here by the exponent ``$0$''). When restricted to a ``quadratic form of M. Hermite,'' in Loewy's terms, that is, when $a_{ik}=a_{ki}^0$, it says that the form can be transformed into itself only by transformations whose characteristic polynomials have simple elementary divisors and roots of modulus 1. Loewy then stated the above theorem, and used it in particular to complete a previous study by Picard of the finite groups of ternary linear transformations (\cite{Picard1887}). For each such group, except one, Picard had displayed a quadratic form left invariant by the group, and Loewy provided an explicit invariant Hermitian form for the remaining case. Following immediately after the publication of his note, there was unleashed a flood of publications.

On August 9, Lazarus Fuchs communicated to the French Academy a note pointing out that what Loewy had published ``\emph{without proof}'' (Fuchs's emphasis) was a special case of his own results presented at the Berlin Academy exactly one month earlier, on July 9 (\cite{Fuchs1896}). Now Fuchs had, for quite some time, been studying differential equations of the type:

\begin{equation*}
\frac{d^n \omega}{dz^n}+q_1\frac{d^{n-1} \omega}{dz^{n-1}} +\cdots + q_n \omega=0
\end{equation*}

\noindent where the coefficients $q_i$ are uniform functions of the variable $z$, with a finite number of poles. Choosing a set of $n$ linearly independent solutions in the neighbourhood of a singularity,  the new solutions one  obtained when the variable $z$ describes a circuit around the singularity can be expressed as the image by a linear transformation of the original ones, and Fuchs, among others, had studied these monodromy transformations, in particular their ``fundamental equations'' (for us, the characteristic equations of the matrices associated to these monodromy transformations) (\cite{Gray1984}). In his July presentation, Fuchs had stated that, under several assumptions (in particular, quite unnecessarily, that the roots of at least one fundamental equation should be distinct), there exists a linear combination of a fundamental system of solutions $\omega_i$ of the differential equation 
\begin{equation*}
\phi = A_1 \omega_1\omega_1'+ A_2 \omega_2\omega_2' +\cdots + A_n \omega_n\omega_n'
\end{equation*}
\noindent (${\omega_i}'$ being here the conjugate function), with determined real coefficients $A_i$, which is unaltered by the group of monodromy. For algebraically integrable differential equations, the group is finite and Fuchs had also used his theorem to complete Picard's work on ternary forms. 

Felix Klein took the opportunity of the annual meeting of the Deutsche Mathe\-ma\-ti\-ker-Vereinigung, from September 21 to September 26,  in Frankfurt  to  present a one-page paper which added another author  and another filiation to the theorem (\cite{Klein1896}). First of all, Klein recalled his own 1875 work (\cite{Klein1875}) where he had explained how to interpret a finite group of complex linear transformations on two variables as a group of real quaternary collineations of the ellipsoid, $x^2+y^2+z^2-w^2=0$; and that this group necessarily fixes a point within the ellipsoid, thus providing a finite group of (real) rotations around a fixed point. This was the basis of his own classification, for the binary case. The ternary case he attributed not only to Picard, but also to Hermann Valentiner. Valentiner, who after a thesis on space curves had gone to work for a Danish insurance company, while still contributing to mathematics, had indeed published in 1889 a book on the classification of finite binary and ternary groups of linear transformations (including the now-called Valentiner group)  (\cite{Valentiner1889}). Then, after a nod to Loewy's note, Klein devoted the remainder of his presentation to another proof that Eliakim Hastings Moore, from Chicago, had communicated to him: For any Hermitian definite form, the sum of its transformations by the (finitely) many elements of the group is still a Hermitian definite form and it is fixed by the group. In the written version of his communication, Klein added that Moore had indeed spoken about his theorem at a mathematical meeting in Chicago on July 10 (with a written version published locally on July 24)!\footnote{Archives of the Math Club, Box 1, folder 2, p.66; see (\cite[399]{PR1994}) for the Chicago environment of this theorem and a slightly different datation. On the argument, see (\cite[511-512]{Hawkins2013});  it is remarkable to us that Klein felt it necessary  to explain that such a procedure would not necessarily work if the group was infinite.} He also alluded to Fuchs's work without more details (given the past tensions between Fuchs and Klein, one may think that this vague recognition was not completely satisfactory to Fuchs (\cite{Gray1984})).

Both Moore and Loewy published an extended version of their respective work in 1898, in the same issue of the \emph{Mathematischen Annalen}, of which Klein was editor-in-chief. Both men analysed the literature, and in particular  underlined Fuchs's  superfluous condition to refute his claim to the theorem (Loewy emphasising that, above all, the definite character of the invariant form was never even alluded to by Fuchs). Still their viewpoints were quite different, as can already been inferred from the proofs themselves and their backgrounds. Moore attributed to ``the analytic phrasing in terms of binary groups of Klein's invariant point'' (in Klein's 1875 paper) his own discovery of the universal invariant Hermitian form (\cite{Moore1898}). After he had proved the existence of this form, he proceeded to deduce from it the theorem that a $n$-linear homogeneous transformation of finite order $p$ can be written with adequate new variables as the multiplication by $p$\textsuperscript{th}-roots of unity, a theorem for which he quoted no less than five other proofs. One of these, included just after his own paper in \emph{Mathematische Annalen}, was due to his Chicago colleague Heinrich Maschke (\cite{Maschke1898}); during this period, Maschke also worked on the theory of groups of linear transformations and he had even delivered a survey lecture at the Chicago Mathematical Club in May 1897. He would soon use Moore's theorem to study quaternary groups of transformations (\cite[396-401]{PR1994}), before extending these results to his now celebrated statement on the representation of finite groups (\cite[512]{Hawkins2013}).

Loewy (who developed his work even further in his 1898 \emph{Habilitation} at the Albert-Ludwigs-Universit\"at in Freiburg) focussed not on the finite groups, but on bilinear forms (\cite{Loewy1898}); his aim was the study of the transformation of a bilinear form with conjugate complex variables (with non-zero determinant) into itself,  in direct continuation of the theme of his thesis. His framework and his main reference is  Georg Frobenius's work, to which he borrowed in particular his symbolic methods (\cite{Frobenius1878}). Such a symbolism, although not linked to a matricial setting, would allow him for instance to characterise a Hermitian form $S$ by the equation $\bar{S}' =S$ (the bar now designates the complex conjugation, the prime the transpose). In the last part of his paper, Loewy also addressed the issue of a reduction theory for Hermitian forms in $n$ variables. It provided him with information on the characteristic polynomial of the linear transformations fixing a given Hermitian form, thus generalising a theorem of \cite{Frobenius1883}.
Although their name was not yet set, the bilinear and the quadratic viewpoints have thus merged and Hermitian forms have become at the end of the nineteenth century a familiar object in the nascent area of linear algebra.

\section{Hermitian forms as geometric  objects}
As seen in Picard's work, Hermitian forms have been connected to a geometric setting already in the 1880s. But it is in the framework defined  by Corrado Segre who, in 1890,  opened new vistas in complex geometry, that they will take a key place as geometric objects (\cite{Segre1890a, Segre1890b}). A starting point for Segre was Karl Von Staudt's project  ``to
make the geometry of position  an autonomous science,  which does not need measure'' (\cite{Staudt1847}). Instead of using cross-ratios, as his predecessors had done, von Staudt based his geometry only on the concept of harmonicity (corresponding to the case where the cross-ratio is -1), which he defined by a purely geometric construction. Two geometric entities whose respective elements are put into correspondence are then said to be projectively related if the  elements of the second entity corresponding to a harmonic set of elements of the first entity are also an harmonic set. A key concept  was that of two entities in involution: for instance, two projectively related ranges of points on the same line are said to be in involution if the corresponding point of a point $A'$ corresponding to a point $A$ is the point $A$ itself; in modern terms, it means of course that the transformation from the first range of points to the second is of order 2, but for von Staudt, to be projectively related or to be in involution were \emph{relations} between  geometric entities ; there were no transformations, as mathematical objects \emph{per se}. For some involutory relations, there exists fixed (real) elements, for others none. For instance, any line  through a point $O$ inside a circle---let us remind that in this projective setting, a ``circle'' is simply a projective line---cuts the circle in two points $A$ and $A'$ ; the ``sheaf'' of all these lines, in von Staudt's terminology, defines an involutory relation between all the points $A$ of the circle and their corresponding points $A'$. This relation has no fixed point. On the other hand, if the point $O$ is outside the circle, the analogous relation has two fixed points, that is, the points of contact with the circle of the two tangents passing through $O$.  In order ``to give to projective geometry the same generality as analysis,'' von Staudt found it necessary to include imaginary (that is, complex) points into geometry, but: ``Where, everybody asks, is the imaginary point when one abstracts oneself from the system of coordinates?'' (\cite{Staudt1847}). His first solution was  to associate  to an involutory situation on a real projective line with no real fixed points \emph{a pair of conjugate imaginary points}; a few years later, he  refined it, distinguishing between the two conjugates by taking into account an additional direction, restricted himself to projective relations leaving the direction invariant and thus extended all his constructions and incidence relations to complex points (\cite{Nabonnand2008}).

Despite the uncompromising austerity of von Staudt's treatise and the complete absence of figures, the reception of this approach was particularly good as a prerequisite for courses in technical drawings and in engineering schools, as advocated in particular by Carl Culman and Theodor Reye at Z\"urich's Polytechnikum (the future ETH) (\cite{Scholz1989}): 

\begin{quote}
But still the engineer, and everybody else who wishes to become familiar with his ideas, must continually exert his power of imagination in order actually to see the object intended to be represented by the lines of a drawing which is not at all intelligible to the uninitiated. […] One principal object of geometrical study appears to me to be the exercise and the development of the power of imagination in the student, and I believe that this object is best attained in the way in which Von Staudt proceeds,
\end{quote}

\noindent wrote Reye at the beginning of his own introduction to projective geometry  (\cite{Reye1866}). The topic gained a particular popularity in  Italy after the \emph{Risorgimento}, where the teaching of geometry was favoured both for national and epistemological reasons, while engineering training responded to the needs of modernization (\cite{Bottazzini1994}). From Luigi Cremona's \emph{Elementi di geometria projettiva} and Antonio Favaro's \emph{Lezioni di statica grafica}  in 1873 
to Achille Sannia's \emph{Lezioni di Geometria Proiettiva} (1891) to Francesco Severi 's \emph{Geometria projettiva} (1904), courses blossomed in Turin, Padua, Bologna, Rome or Palermo (\cite{CR2014}). In Turin, Corrado Segre, besides his own lectures on projective geometry,\footnote{His lecture notes are available on the beautiful website devoted to him by Livia Giaccardi, \url{http://www.corradosegre.unito.it/I21_30.php} see in particular Quadreno 22.}   significantly organised the translation into Italian of both von Staudt’s \emph{Geometry of position} and Klein’s Erlangen program, by Mario Pieri (1889) and Gino Fano (1890) respectively. 
A decisive change with respect to von Staudt's work is that, for Segre (as for Klein),  a projectivity is a transformation (e.g. a collineation, a homography), a mathematical object, and not only a way of relating  entities. In complex geometry, projectivities preserve harmonic relations, but the reciprocal is false, for instance the complex conjugation is not a projectivity. While von Staudt had chosen to restrict the relations he considered, Segre introduced what he called anti-projectivities, as well as their 2- or 3-dimensional companions (anticollineations, antipolarities,…) and studied projectivities and anti-projectivities alike. Analytically, a projectivity  in dimension 1 corresponds to a homography $$z \rightarrow z'= \frac{\alpha z + \beta}{\gamma z + \delta}$$

\noindent while an anti-projectivity corresponds to : $$z \rightarrow z'= \frac{\alpha \bar{z} + \beta}{\gamma \bar{z} + \delta}$$

\noindent with $\alpha, \beta, \gamma, \delta$ complex numbers such that $\alpha\delta-\beta \gamma \not= 0$ (the bar indicates the complex conjugation).

It is then easy to see that, if one puts  $z=\frac{x}{y}$, the equation giving the fixed points of such an antiprojectivity is of the type:

$$F(x,y)= ax\bar{x}+bx\bar{y}+\bar{b} \bar{x}y+c y\bar{y}=0$$

\noindent Similarly,  antipolarities are associated to  forms with 4 variables.  One recognises the forms that, in Segre's words, ``have also already been introduced in number theory thanks to M. Hermite, M. Picard and others,''  that is, again, Hermitian forms  (\cite[275]{Brigaglia2016}). Segre explained in a letter to Adolf Hurwitz in June 1894: ``I study there, among the other hyperalgebraic entities, those that I call hyperconic, hyperquadric [i.e. the locus of $F(x_1, x_2, x_3, x_4)= 0$], etc. which are analytically represented by the equations of Hermite; and thus they give the geometric equivalent of the forms of Hermite'' (\cite{Brigaglia 2016}). 

According to Klein's Erlangen program, the group of projective and antiprojective transformations which fix, say, an hyperquadric, defines a geometry, and it is also
possible to apply Cayley's idea, as extended by Klein, of using a (hyper)quadric to define a metric on this projective complex space. This was done almost simultaneously and independently by both Guido Fubini and Eduard Study.  Fubini had been a student of Ulisse Dini and Luigi Bianchi at Pisa and Bianchi included results from Fubini's thesis on Clifford parallelism in his celebrated 1902 book \emph{Lezioni di geometria differenziale}. A year later, Fubini, then in Catania, turned to Hermitian forms  (\cite{Fubini1982}). Starting from the interpretation of the group of linear transformations leaving invariant a quadratic form as a group of space motions, he was seeking an analogous construction for Hermitian forms. Fubini considered an Hermitian form algebraically equivalent to $x_1 x_1^0+ x_2 x_2^0+\cdots +x_{n-1} x_{n-1}^0- x_n x_n^0$ (here again the exponent $0$ indicates the complex conjugation) and the associated group transforming a ``hypervariety'' of the type $\sum_1^{n-1} u_i u_i^0 - 1=0$ into itself. For two points, $u$ and $\bar{u}$ (the bar here does \emph{not} designate the complex conjugation), he introduced the quantity

\begin{equation}\label{fubini}
R_{u\bar{u}} =\frac{(\sum_1^{n-1} u_i \bar{u}_1^0-1)(\sum_1^{n-1} u_i ^0\bar{u}_1-1)}{(\sum_1^{n-1} u_i {u}_1^0-1)(\sum_1^{n-1} \bar{u}_i ^0\bar{u}_1-1)} -1.
\end{equation}

This is invariant under the group of transformations, real and equal to $0$ within the hypervariety only when the two points coincide; it was thus legitimate to call it a (pseudo)-distance between the two points. \cite{Fubini1903} used it primarily with arithmetical and analytical applications in view, but he then also refined his construction in order to interpret the Hermitian form as a metric for a complex space (\cite{Fubini1904}).

While Fubini refered mostly to Picard, Poincar\'e and Klein-Fricke, Study's idea  was to develop this study within the framework defined by Segre. Study had worked before on invariant theory and quaternions and had just published in 1903 a book on \emph{Geometrie der Dynamen}, using biquaternions and geometrical tools to study mechanical forces (\cite{Hartwich2005}). He had an extensive program about the complex realm and in an article published in \emph{Mathematische Annalen} (\cite{Study1905}), he also defined Hermitian metrics and distances. 

From a ternary indefinite Hermitian form, $(x\overline{x})=x_1\overline{x_1}-x_2\overline{x_2}-x_3\overline{x_3}$ for instance (here the bar does designate the complex conjugation), Study, following Segre's concepts if not his terminology, defined a ``Hermitian point-complex'' by $(x\overline{x})=0$ (this is Segre's \emph{iperconica}) and the inside of the point-complex by the condition $(x\overline{x})>0$ (by the very definition of a Hermitian form, its values are real). He was then able to define a hyperbolic Hermitian metric and the (real) distance of two points inside  the point-complex. Under an adequate normalization,  the distance between two points $x$ and $y$ is 

\begin{equation*}
(x,y)= 2 \cosh^{-1} \frac{\sqrt{(x\overline{y})(\overline{x}y)}}{\sqrt{(x\overline{x})}\sqrt{(y\overline{y})}},
\end{equation*}

\noindent Study showing then that the distance between two points is the length of the geodesics linking them. He also developed  the case of  an elliptic Hermitian metric, based this time on a definite Hermitian form.
This setting would then be developed by Julian Coolidge, Wilhelm Blaschke and of course Erich K\"ahler, Jan Schouten and \'Elie Cartan in the 1920s and 1930s.

\section{Open questions: back to long-term histories}
As we have seen, while in the 1850s, what will be called Hermitian forms were a subcategory of quadratic forms, they became later a category \emph{parallel} to that of quadratic forms—and still later the encompassing category of those quadratic forms which are a particular case (those with complex conjugation reducing to identity). Alternatively, it may be tempting to summarise the history of Hermitian forms, of which we have just sketched some key episodes, as a thread from number-theoretical objects to algebraic ones to geometric ones. However, one would obviously lose several components we have met — be they pedagogical settings, or epistemological views, putting observations and classifications together at the core of mathematical activities, or disciplinary structures. In the configuration which allows Hermitian forms to appear as geometric objects, there appears the renewal of Italian mathematics after the unification of the country, a new program for a geometry as general as analysis, a concern for the representation of complex numbers by real objects, and an emphasis on transformations. If close reading is required to display how these components were articulated, the question is still open of their extent, both in time and in scope, and of their precise role in the dynamics of the changes.

The difficulty is hardly reduced if one focuses on components which are obviously relevant, for instance here, the issue of geometrisations: several types of geometries occurred (projective, non Euclidean, differential) and they played several roles, from supporting intuition to providing a concrete existence to groups. As Fubini (\cite[2]{Fubini1904}) explained, for instance, ``the introduction of [his] new metric is something totally different from a simply formal thing. It allows us to have recourse to the intuition and the procedures of geometry to solve an algebraic question''. Moreover, geometry did not go unchallenged: Valentiner remarked at the beginning of his book on finite groups of transformations that while Klein had already handled the binary case, he, Valentiner, was going to do it again, but this time purely algebraically, in order to free it from geometrical considerations (\cite[5]{Valentiner1889}). The relation between geometry and classification programmes was established in Segre’s case, according to his former student and colleague Alessandro (\cite{Terracini1926}); it was however still an open question whether to connect it or not with more general views on mathematics as a science, as was the case for Hermite or Cayley (\cite{LP2016}).

Besides these questions about specific components, an important range of open questions concerns the links which may be used to relate the various episodes and components, in particular those which can help us to understand how information or results were circulating. An obvious type of link, which has been used to delineate preceding configurations, are cross-references among articles. They indicate in particular the role played by Picard in the 1880s to put Hermitian forms back into the centre of the stage. One might also think of Cayley’s definition of a metric in a projective space, linked to his work on invariants in the 1860s, and which will be, in a form mainly mediated by Klein, a model for the formulas used in the complex case at the end of the century. 

Conflicts of priority often gave rise to such cross-references, but they were, at least in our case, particularly ephemeral. Issues were very soon settled in the literature: there will be a ``Fubini-Study'' metric, the theorem on the invariance of a Hermitian form by a finite group of linear transformations will be attributed to all the authors we have discussed above (\cite[341]{Meyer1899}, \cite[209]{MBD1916}). Yet locating and studying such configurations remains of interest, as they help us reveal more decisive aspects: such as, for instance, the confluence of different approaches set in motion to describe the finite subgroups of $GL_n (\mathbb{C})$, or the problems raised by the extension of geometry to the complex realm, or the role of Klein as a intermediary between Göttingen, Italy and the United States (\cite{PR1994, CCGGMV2016}). That persons from different branches of mathematics could agree was even interpreted by some as a sign of modernity. Francesco Brioschi commented for instance in 1889: ``The characteristic note of modern progress in mathematical studies can be recognised in the contribution that each special theory—that of functions, of substitutions, of forms, of transcendents, geometrical theories and so on---brings to the study of problems, where in other times only one seemed necessary'' (\cite[375]{Archibald2011}).

However, transmissions may also operate without explicit bibliographical indications, through correspondence or, as we have seen, through communications at the meetings of a mathematical society or via publication in a given journal (the role of \emph{Mathematische Annalen} is here particularly interesting). Certain links functioned on a much larger scale and to pinpoint them concretely is more delicate: how a new object in mathematics is validated or how the new ideas of linear algebra irrigated a variety of topics during the second half of the nineteenth century are but two examples.

Further, breaks or discontinuities should also be studied more closely. This is particularly important because most theories in the dynamics of history of science rely on the issue of breaks (from Thomas Kuhn’s scientific revolutions to Imre Lakatos’s research programmes) and do not seem to be appropriate in the case of mathematics. It may happen that different aspects change simultaneously and there have been some attempts to speak in this case of a ``revolution'' (\cite{Gillies1992}). Very often however the changes are not simultaneous (\cite{GilGui2015}). While indifference is of course the most obvious cause of decline, rediscovery or recycling into another theory are quite frequent in modern times and the circumstances of such a rebirth are often puzzling. We have only a few studies of such phenomena: that of invariant theory, dead (\cite{Fisher1966}), but then born again, ``like an Arabian phoenix arising from its ashes’' (\cite{Rota2001}). Or that of the ``theory of order’' that Louis Poinsot, inspired by the relations among roots of unity, promoted at the beginning of the century (\cite{Boucard2011}), and which disappeared and reappeared several times (still attached to Poinsot's name), in the theory of equations, but also in ornamental architecture (\cite{BouEck2015}), intersecting the long-term changes in the theory of tactics, a field which mixed what we would now describe as combinatorics and group theory (\cite{Ehrhardt2015})). In examples ranging from Descartes’s curves to Fourier series, Alain \cite {Herreman2013} analysed how concepts which had already appeared informally in mathematics receive, in what he calls ``inaugural texts,’’ a formal definition and a name, thus acquiring, so to speak, mathematical citizenship; but such a general study of a specific type of discontinuity, operating here in particular at a semiotic level, remains exceptional.

It may sometimes be necessary to change one’s scale of observation in order to understand a discontinuity. A striking testimony is given by Max Born's 1914 study of the diamond (\cite{Born1914}). Born used a determinant with complex elements, such that two elements symmetrical with respect to the main diagonal are complex conjugates, commented that the determinant is ``real despite its complex appearance’’ and acknowledged in a footnote that it thus derived from a Hermitian form. As well-known, Born would remember this again a decade later, offering to Hermitian forms and matrices a spectacular role in the new quantum theory (\cite{BJ1925}). In the 1850s, the program of classification and reduction of equations and forms had generated various concepts: characteristic polynomials and their roots, which we have met several times, linear transformations or matrices (\cite{Hawkins1977, Hawkins1986, Brechenmacher2010, Hawkins2013}). That the ``secular equation’' which describes the long-term perturbations of planetary motion is our characteristic polynomial for a symmetric matrix and has only real roots had indeed interesting physical applications; it was an incentive to work on the complex case for \cite{Clebsch1860} and \cite{Christoffel1864}, whom I have mentioned above. Secular equations again appeared in a physical environment in Born’s paper, but without explicit link to these previous works. The connections between Born and Jacob Rosanes, who studied the transformation of quadratic forms in the 1870s and 1880s, are on the other hand well attested (\cite[119]{MR1982}) and Born also quoted more recent algebraic papers: this suggests a wider, less specific, relevant algebraic environment in which matrices and forms were now available independently of the narrow topic of the secular equation or of its other avatars.

In general history, following Fernand Braudel's suggestion to take into account not only the short-term events of political life, but also phenomena deployed over a more extended period of time, be they economic trends or even geological transformations, the ``long term'' has come to mean a time of immobility with respect to human actions. This is not the case in mathematics. A long term thread, based on a simple retrospective identification of concepts, or results, would be a snare. There is no easy line from Hermite to Picard to Moore and Study to Born. Long term histories thus require us to study not only local configurations but also the various ways in which they are, or not, connected. The links we have mentioned, whatever be their duration in time or the scale at which they operate, are constructed through the ``art’’ of mathematicians; and this may be true also of the discontinuities. Insofar as they involve the work of mathematicians, they are---or should be---a primary concern for the history of mathematics.

\nocite{*}
\bibliography{bibgoldsteindef}

\end{document}